\documentclass[12pt,leqno]{article}
\usepackage{amsthm,amsfonts,epsfig}
 \usepackage{latexsym, amssymb}

\def\ds{\displaystyle}\def\scs{\scriptstyle}

\def\ni{\noindent}\def\ov{\over}\def\p{\partial}
\def\smoc{{\cal C}^\infty_c}
\def\td{\tilde}\def\la{\langle}\def\ra{\rangle}
\def\ep{\epsilon}\def\om{\omega}\def\lam{\lambda}
\def\Gam{\Gamma}\def\gam{\gamma}
\def\ovl{\overline}

\newtheorem{theorem}{Theorem}[section]

\newtheorem{lemma}{Lemma}[section]

\title{Bessel Integrals and Fundamental Solutions for a Generalized
      Tricomi Operator} 
\author{J.~Barros-Neto\\
Rutgers University, Hill Center\\
110 Frelinghuysen Rd, Piscataway, NJ 08854-8019\\
e-mail: jbn@math.rutgers.edu
\and
Fernando Cardoso\footnote{Partially supported by CNPq (Brazil)}\\
Departmento de Matem\'atica, Universidade Federal de Pernambuco\\
50540-740 Recife, Pe, Brazil\\
e-mail: fernando@dmat.ufpe.br} 
\date{} 
\begin{document}  
\maketitle
\begin{abstract}
The method of partial Fourier transform is used to find explicit
formulas for two remarkable fundamental solutions for a generalized
Tricomi operator. These fundamental solutions reflect clearly the
mixed type of the Tricomi operator. In proving these results, we
establish explicit formulas for Fourier transforms of some functions
involving Bessel functions.
\end{abstract}

\section{Introduction}\setcounter{equation}{0}\label{in}
Consider the generalized Tricomi operator
\begin{equation}\label{eq1}
P=y\Delta +{\p^2\over\p y^2},
\end{equation}
where $\Delta=\sum_{j=1}^n\ds{\p^2\ov\p x_j^2}$ is the Laplace
operator. Our aim is to find, via partial Fourier transform with
respect to $x=(x_1,\ldots,x_n)$, fundamental solutions relative to an
arbitrary point $(a,0)$ on the hyperplane $y=0$ in ${\mathbb
R}^{n+1},$ that is, distributions that are solutions to
$$
Pu=\delta(x-a,y),
$$ 
where $\delta(x-a,y)$ is the Dirac measure concentrated at $(a,0).$
Since $P$ is invariant under translations parallel to that hyperplane,
it suffices to consider the case when the Dirac measure is
concentrated at the origin.

If $n=1,$ then  (\ref{eq1}) is the classical Tricomi operator,
\begin{equation}\label{eq2} 
{\cal T}=y{\p^2\over\p x^2}+{\p^2\over\p y^2}.
\end{equation}
It is known that, for this operator, the equation $9x^2+4y^3=0$
defines the two characteristic curves that originate at the
origin. They divide the plane in two disjoint regions, namely,
$$
D_+=\{(x,y)\in{\mathbb R}^2 : 9x^2+4y^3>0\},
$$
the region ``outside'' the characteristics and 
$$
D_-=\{(x,y)\in{\mathbb R}^2 : 9x^2+4y^3<0\},
$$
the region ``inside'' the characteristics. Note that $D_-$ is entirely
contained in the hyperbolic region $y<0.$

In the paper \cite{bg} it was shown the existence of the following two
fundamental solutions relative to the origin:

\begin{equation}\label{eq2a}
\hskip .5cm F_+(x,y) =\left\{\begin{array}{ll}
         -\ds{\Gam(1/6)\ov3\cdot2^{2/3}\pi^{1/2}\Gam(2/3)}(9x^2+4y^3)^{-1/6}
         & \mbox{in $D_+$}\\
\\
          \ \ \  0 & \mbox{elsewhere}
\end{array} \right .
\end{equation}
and 
\begin{equation}\label{eq2b}
F_-(x,y) = \left\{\begin{array}{ll}
         \ds{3\Gam(4/3)\ov2^{2/3}\pi^{1/2}\Gam(5/6)}|9x^2+4y^3|^{-1/6}
         & \mbox{in $D_-$}\\ \\ 0 & \mbox{elsewhere}.
\end{array} \right .
\end{equation}

The support of $F_-,$ the closure of $D_-,$ is, except for the origin, entirely contained in the
hyperbolic region ($y<0$), while the support of $F_+,$ the closure of $D_+,$ consists of the whole elliptic
region ($y>0$), the parabolic region ($x$-axis), and extends to the hyperbolic region up to and included the
characteristic curves. The method used in \cite{bg} to prove these results was based upon the property that
solutions to the equation ${\cal T}u=0$ are invariant under the dilation $d_t(x,y)=(t^3x,t^2y)$ in ${\mathbb
R}^2.$ 

In \cite{bn} we proved formulas (\ref{eq2a}) and (\ref{eq2b}) by using
partial Fourier transform with respect to the $x$ variable.  In this
paper we extend the results of \cite{bn} to the generalized Tricomi
operator (\ref{eq1}). Since the dimension of the space variable $x$ is
now $n>1,$ technical difficulties in evaluating Fourier transforms
involving Bessel functions do occur. We circumvent them by calculating
integrals of the type
$$
I_\ep(a,b)=\int_0^\infty e^{-\ep t}t^{-\lam}J_\mu(at)J_\nu(bt)\, dt
$$
(Section \ref{ft}). As a result, we obtain for the operator
(\ref{eq1}) the following fundamental solutions:
$$
F_+(x,y) = \left\{\begin{array}{ll}
\ds-{3^{n-2}\ov
2^{2/3}\pi^{n/2}}{\Gamma({n\ov2}-{1\ov3})\ov\Gamma(2/3)}
      \,(9|x|^2+4y^3)^{{1\ov3}-{n\ov2}}
       & \mbox{in $D_+^n$}\\
\\
   \hskip 1cm 0  &\mbox{elsewhere}, 
\end{array} \right.
$$
whose support is the closure of the region $D_+^n=\{(x,y)\in{\mathbb
R}^{n+1}: 9|x|^2+4y^3>0\}$ and
$$
F_-(x,y) = \left\{\begin{array}{ll}
\ds{3^n\Gamma({4\ov3})\ov2^{2/3}\pi^{n/2}\Gamma({4\ov3}-{n\ov2})}
      \,|9|x|^2+4y^3|^{{1\ov3}-{n\ov2}}
       & \mbox{in $D_-^n$}\\
\\
 \hskip 1cm 0 & \mbox{elsewhere} 
\end{array} \right.
$$
supported by the closure of the region $D_-^n=\{(x,y)\in{\mathbb
R}^{n+1}: 9|x|^2+4y^3<0\}.$

These fundamental solutions clearly generalize formulas (\ref{eq2a})
and (\ref{eq2b}). In addition, we also show that no matter the parity
of the dimension $n,$ $F_-(x,y)$ {\bf does not} have support on the
boundary of the region $D_-$ in the hyperbolic region and hence, in
contrast with what happens for strictly hyperbolic operators -- like
the wave operator -- the Huyghens principle does not hold for the
generalized Tricomi operator.

It is well known \cite{tv} that, for the wave operator, the regularity
of the fundamental solutions degenerates as the dimension $n$
increases: a locally constant function, when $n=1;$ an absolutely
continuous measure relative to the Lebesgue measure, when $n=2;$ a
measure carried by the surface of the forward light-cone, when $n=3;$
and so on. However, for the generalized Tricomi operator they always
remain locally integrable.

The plan of this paper is as follows. In Section \ref{pr} we reduce,
via partial Fourier transform, the original problem to an equivalent
one of finding fundamental solutions for a second order ordinary
differential equation and show how these can be represented in terms
of Airy functions or Bessel functions.  In Section \ref{ft} we obtain
explicit formulas for Fourier transforms of the functions
$|\xi|^{\pm\nu}J_\nu(|\xi|),$ $|\xi|^{\nu}K_\nu(|\xi|),$ and
$|\xi|^{\nu}N_\nu(|\xi|),$ $\xi=(\xi_1,\ldots,\xi_n),$ $n\geq1.$ These
formulas are used to construct the fundamental solution supported by
$\ovl {D_-^n}$ (Section \ref{hy}) and the fundamental solution
supported by $\ovl {D_+^n}$ (Section \ref{el}). In the Appendix the
reader will find the definitions of the Bessel functions $J_\nu(z),$
$I_\nu(z),$ $K_\nu(z),$ and $N_\nu(z),$ and a list of properties that
are needed throughout this work. We also recall the definition of
hypergeometric functions and some of their main properties.

The method of finding solutions, via partial Fourier transform, to the
equation
\begin{equation}\label{eq1a}
Pu=y\Delta u+{\p^2 u\over\p y^2}=f
\end{equation}
with, say, $f\in\smoc({\mathbb R}^{n+1})$ is a natural one although
not particularly new. In the monograph \cite{gro},
R.~J.~P.~Groothuizen exhibits a solution $u$ in terms of a Fourier
integral operator whose symbol is obtained by the partial Fourier
transform analysis employed in this paper. He also obtains fundamental
solutions given by Fourier integral operators. Our results are more
precise since we give explicit formulas for the fundamental solutions
and these are tempered distributions.  Consequently, we can obtain more
transparent representations for the solution $u$ to the equation
(\ref{eq1a}) as a convolution of $f$ with any of the fundamental
solutions here described.

\section{Preliminaries}\setcounter{equation}{0}\label{pr}

Consider the more general problem of finding fundamental solutions for
the operator (\ref{eq1}) relative to an arbitrary point $(0,b)$ on the
$y$-axis. That is, one wishes to find distributional solutions to the
equation
\begin{equation}\label{eq3}
y\Delta F+{\p^2\over\p y^2}F=\delta(x)\otimes\delta(y-b).
\end{equation}
Partial Fourier transform with respect to $x$ reduces this problem
into finding fundamental solutions to the ordinary differential
equation
\begin{equation}\label{eq4}
\td F_{yy}-y|\xi|^2\td F=\delta(y-b),
\end{equation}
where $\td F$ denotes the partial Fourier transform of $F,$
$|\xi|^2=\sum_{j=1}^n\xi_j^2,$ and $\delta(y-b)$ is the Dirac measure
concentrated at $b.$

One way of solving (\ref{eq4}) (see \cite{tv}) consists in selecting
two linearly independent solutions $U_1(\xi,y)$ and $U_2(\xi,y)$ to
the homogeneous equation
\begin{equation}\label{eq4a}
\td u_{yy}-y|\xi|^2\td u=0
\end{equation}
so that their Wronskian at $y=b$ (which in the case under
consideration is the same as the Wronskian at $y=0$) is equal to $-1$
and in defining
\begin{equation}\label{eq4b}
\td F(\xi,y;b) = \left\{\begin{array}{ll}
         U_2(\xi,b)U_1(\xi,y) & \mbox{if $y\geq b$}\\
          \\
         U_1(\xi,b)U_2(\xi,y) & \mbox{if $y\leq b$}.
\end{array} \right.
\end{equation}
It is a matter of verification that $\td F(\xi,y;b)$ is a solution to
(\ref{eq4}).

Equivalently \cite{gs}, one can obtain a fundamental solution to the
equation (\ref{eq4}) by finding two linearly independent solutions to
the homogeneous equation (\ref{eq4a}): $u_1(\xi,y;b),$ defined for
$y>b,$ and $u_2(\xi,y;b),$ defined for $y<b,$ so that

(i) the limit of $u_1$ as $y\to b+$ equals the limit of $u_2$ as $y\to
b-$
\medskip

\ni and

\medskip
(ii) the limit of $u_{1,y}$ as $y\to b+$ minus the limit of $u_{2,y}$
as $y\to b-$ is equal to 1. The function $\td F(\xi,y;b)$ is now
defined by
\begin{equation}\label{eq4c}
\td F(\xi,y;b) = \left\{\begin{array}{ll}
         u_1(\xi,y;b) & \mbox{if $y\geq b$}\\
          \\
         u_2(\xi,y;b) & \mbox{if $y\leq b$}.
\end{array} \right.
\end{equation}
In what follows we will use at our convenience either one of these two
expresions for $\td F(\xi,y;b).$

Returning to the original problem (\ref{eq3}), if we can choose
$U_1(\xi,y)$ and $U_2(\xi,y),$ in formula (\ref{eq4b}), or
$u_1(\xi,y;b)$ and $u_2(\xi,y;b),$ in formula (\ref{eq4c}), so that
$\td F(\xi,y;b)$ is a tempered distribution with respect to
$\xi=(\xi_1,\ldots,\xi_n),$ then its inverse Fourier transform of $\td
F(\xi,y;b)$ defines a fundamental solution to the operator
(\ref{eq1}).

We now proceed to find linearly independent solutions to the
homogeneous equation (\ref{eq4a}).  The change of variables
$z=|\xi|^{2/3}y$ transforms (\ref{eq4a}) into the classical Airy's
equation $u''-zu=0.$ Two linearly independent solutions to that
equation are $Ai(z)$ and $Bi(z)$ respectively called {\em the Airy
functions of the first and second kinds.} It is known \cite{lb} that
these two functions can be represented in terms of Bessel functions of
order $\pm1/3$ as follows. If $|\arg z|<(2\pi/3),$ then

\begin{equation}\label{eq5}
\hskip 1cm \left\{\begin{array}{ll}
Ai(z)={z^{1/2}\over 3}[I_{-1/3}({2\over3}z^{3/2})-I_{1/3}({2\over3}z^{3/2})]
     = {1\over\pi}({z\over 3})^{1/2}K_{1/3}({2\over3}z^{3/2}),\\
\\
Bi(z)=({z\over3})^{1/2}[I_{-1/3}({2\over3}z^{3/2})+I_{1/3}({2\over3}z^{3/2})].
\end{array} \right. 
\end{equation}
If $|\arg z|<(2\pi/3),$ then
\begin{equation}\label{eq6}
\left\{\begin{array}{ll}
Ai(-z)={z^{1/2}\over3}
[J_{-1/3}({2\over3}z^{3/2})+J_{1/3}({2\over3}z^{3/2})],\\
\\
Bi(-z)=({z\over3})^{1/2}
[J_{-1/3}({2\over3}z^{3/2})-J_{1/3}({2\over3}z^{3/2})].
\end{array} \right.
\end{equation}

In the Appendix I, where a brief review of Bessel functions is
presented, we show that the following relations hold:
\begin{equation}\label{eq6a}
Ai(0)={3^{-2/3}\over\Gamma(2/3)},\qquad Ai'(0)=-{3^{-4/3}\over\Gamma(4/3)}
\end{equation}
and
\begin{equation}\label{eq6b}
Bi(0)={3^{-1/6}\over\Gamma(2/3)},\qquad Bi'(0)={3^{-5/6}\over\Gamma(4/3)}.
\end{equation}
As a consequence, the Wronskian of $Ai(z)$ and $Bi(z)$ evaluated at $z=0$ is
\begin{equation}\label{eq6c}
W(Ai(z),Bi(z))_{|_{z=0}}=1/\pi.
\end{equation}
Indeed, we have
$$
W(Ai(z),Bi(z))_{|_{z=0}}=
\left|\begin{array}{cc}
\ds{3^{-2/3}\over\Gamma(2/3)}&\ds{3^{-1/6}\over\Gamma(2/3)}\\
      \\
\ds-{3^{-4/3}\over\Gamma(4/3)}&\ds{3^{-5/6}\over\Gamma(4/3)}
\end{array}\right|
=\ds{2\cdot3^{-3/2}\over\Gamma(2/3)\Gamma(4/3)}={1\over\pi},\\
$$
because 
\begin{equation}\label{eq6d}
\Gamma(2/3)\Gamma(4/3)={2\pi\over3^{3/2}}.
\end{equation}

We now choose the following two linearly independent solutions to the
equation (\ref{eq4a}):
$$
U_1(\xi,y)=\sqrt\pi|\xi|^{-1/3}Ai(|\xi|^{2/3}y)\quad\mbox{and}\quad
U_2(\xi,y)=-\sqrt\pi|\xi|^{-1/3}Bi(|\xi|^{2/3}y),
$$
and note that, by virtue of (\ref{eq6c}), the Wronskian of
$U_1(\xi,y)$ and $U_2(\xi,y)$ is equal to $-1.$ Next, according to
(\ref{eq4b}), the distribution
\begin{equation}\label{eq7}
\td F(\xi,y;b) = \left\{\begin{array}{ll}
     -\pi|\xi|^{-2/3}Bi(|\xi|^{2/3}b)Ai(|\xi|^{2/3}y) & \mbox{if $y\geq b$}\\
          \\
     -\pi|\xi|^{-2/3}Ai(|\xi|^{2/3}b)Bi(|\xi|^{2/3}y) & \mbox{if $y\leq b$}
\end{array} \right.
\end{equation}
is a fundamental solution to the ordinary differential equation
(\ref{eq4}). A fundamental solution $F(x,y;b)$ to the generalized
Tricomi operator (\ref{eq1}) and relative to the point $(0,b)$ would
then be the inverse Fourier transform of $\td F(\xi,y;b)$ whenever
that Fourier transform exists.

We do not know how to obtain an explicit formula (or formulas) for the
inverse Fourier transform of $\td F(\xi,y;b)$ when $b\neq 0,$ a
problem that merits to be investigated. We conjecture that when $b<0,$
that is, the point $(0,b)$ is in the hyperbolic region, there exists
two fundamental solutions that converge, as $b\to 0,$ to the
fundamental solutions $F_+(x,y)$ and $F_-(x,y)$ defined, respectively,
by the formulas (\ref{eq14}) and (\ref{eq11}). The two fundamental
solutions described in the monograph \cite{gro} do not seem to satisfy
these requirements.

However, when $b=0,$ we will show in the following sections how to
obtain from formula (\ref{eq7}), as well as formulas similar to it,
explicit expressions for fundamental solution to (\ref{eq1}).

\medskip
{\bf The case $n=1$.} For sake of completeness and in order to
motivate our work in the forthcoming sections we briefly present the
results of the paper \cite{bn}, for the Tricomi operator (\ref{eq2}),
that is the case when $n=1.$ If $b=0$ and we take into account the
values of $Ai(0)$ and $Bi(0)$ as given by (\ref{eq6a}) and
(\ref{eq6b}), then formula (\ref{eq7}) simplifies as follows:
\begin{equation}\label{eq8}
\td F(\xi,y) = \left\{\begin{array}{ll}
\displaystyle{-{\pi|\xi|^{-2/3}\over 3^{1/6}\Gamma(2/3)}}Ai(|\xi|^{2/3}y) 
       & \mbox{if $y\geq 0$}\\
          \\
\displaystyle{-{\pi|\xi|^{-2/3}\over 3^{2/3}\Gamma(2/3)}}Bi(|\xi|^{2/3}y)
       & \mbox{if $y\leq 0$},
\end{array} \right.
\end{equation}
where, for simplicity, we wrote $\td F(\xi,y)$ for $\td F(\xi,y;0).$
Now the inverse Fourier transform of both expressions on the
right-hand side of (\ref{eq8}) can be explicitly calculated. To see
this, first consider the representations of $Ai(z)$ and $Bi(z)$ in
terms of Bessel functions as given by formulas (\ref{eq5}) and
(\ref{eq6}):
$$
Ai(|\xi|^{2/3}y)={1\over\pi}\left({|\xi|^{2/3}y\over3}\right)^{1/2}
   K_{1/3}({2\over3}|\xi|y^{3/2}) 
$$
and
$$
Bi(|\xi|^{2/3}y)=\left({|\xi|^{2/3}(-y)\over3}\right)^{1/2}
  [J_{-1/3}({2\over 3}|\xi|(-y)^{3/2})-J_{1/3}({2\over 3}|\xi|(-y)^{3/2})].
$$
Next by introducing the change of variables $s=(2/3)y^{3/2},$ whenever
$y\geq0,$ and $t=(2/3)(-y)^{3/2},$ whenever $y\leq0,$ rewrite formula
(\ref{eq8}) as
\begin{equation}\label{eq9}
\td F(\xi,y) = \left\{\begin{array}{ll}
\ds\alpha\cdot({s\ov|\xi|})^{1/3}K_{1/3}(s|\xi|) 
       & \mbox{if $y\geq 0$} \\
  \\
\ds\beta\cdot({t\ov|\xi|})^{1/3}[J_{-1/3}(t|\xi|)-J_{1/3}(t|\xi|)]
       & \mbox{if $y\leq 0$}
\end{array} \right.
\end{equation}
where $\alpha$ and $\beta$ are constants given by

\begin{equation}\label{eq9a}
\alpha=-{1\over 2^{1/3}3^{1/3}\Gamma(2/3)}\qquad\mbox{and}\qquad
\beta=-{\pi\over 2^{1/3}3^{5/6}\Gamma(2/3)}.
\end{equation}

Theorem 3.1 of \cite{bn} proves that $F(x,y),$ the inverse Fourier
transform of $\td F(\xi,y),$ is then
$$
F={3\over 2} F_+ -{1\over2} F_-,
$$ 
a linear combination of the two fundamental solutions $F_+$ and $F_-$
respectively defined by formulas (\ref{eq2a}) and (\ref{eq2b}). To
prove this theorem, one relies on known formulas (see \cite{gs},
\cite{rg}) for the inverse Fourier transforms of the functions
$|\xi|^{-1/3}J_{1/3}(|\xi|),$ $|\xi|^{-1/3}J_{-1/3}(|\xi|),$ and
$|\xi|^{-1/3}K_{1/3}(|\xi|),$ formulas that need be generalized to the
case $n>1.$

\section{Fourier transforms involving Bessel functions}
  \setcounter{equation}{0}\label{ft} 
As we have indicated in the previous section, we are going to need
explicit formulas for the inverse Fourier transforms of the following
functions: $|\xi|^{\pm\nu}J_\nu(|\xi|),$ $|\xi|^{\nu}K_\nu(|\xi|),$
and $|\xi|^{\nu}N_\nu(|\xi|),$ where $\xi=(\xi_1,\ldots,\xi_n),$
$n\geq1.$ Since, when $n=1,$ explicit formulas for the Fourier
transforms are known and can be found, for example, in \cite{gs} and
\cite{rg}, we concentrate on the case $n>1.$

\bigskip
\ni{\bf 1. The inverse Fourier transform of $|\xi|^{\nu}K_\nu(|\xi|).$ }

Assume that $|\mbox{Re}\ \nu|<1/2$ and define
\begin{equation}\label{k1}
{\cal F}^{-1}[|\xi|^\nu K_\nu(|\xi|)](x)={1\ov(2\pi)^n}
  \int_{\mathbb R_n}e^{i\la x,\xi\ra}|\xi|^\nu K_\nu(|\xi|)\, d\xi.
\end{equation}
The assumption on the real part of $\nu$ secures convergence of the
integral at the origin. On the other hand, in view of the asymptotic
behavior of $K_\nu(|\xi|)$ for large values of $|\xi|,$ the integral
converges absolutely. By introducing spherical coordinates, rewrite
the integral on the right-hand of (\ref{k1}) as
$$
\int_0^\infty r^{n+\nu-1}K_\nu(r)\{\int_{S_{n-1}}e^{i\la rx,\om\ra}\,
d\om\}\,dr.
$$
Since
\begin{equation}\label{k2}
\int_{S_{n-1}}e^{i\la rx,\om\ra}\, d\om=
  {(2\pi)^{\scs n\ov2}\ov|rx|^{{\scs n\ov2}-1}}J_{{\scs n\ov2}-1}(r|x|),
\end{equation}
it follows that
\begin{equation}\label{k3}
\hskip 0.75cm{\cal F}^{-1}[|\xi|^{\nu}K_{\nu}(|\xi|)](x)=
  {|x|^{1-{\scs n\ov2}}\ov(2\pi)^{\scs n\ov2}}
  \int_0^\infty r^{{\scs n\ov2}+\nu}J_{{\scs n\ov2}-1}(r|x|)K_\nu(r)\, dr.
\end{equation}

We now quote the following result found in Watson's treatise \cite{wat}:

\begin{lemma}\label{l0}
If $\mbox{Re}(\mu+1)>|\mbox{Re}\ \nu|$ and $\mbox{Re}\ b>|\mbox{Im}\ a|$ then
\begin{equation}\label{k4}
\int_0^\infty t^{\mu+\nu+1}J_\mu(at)K_\nu(bt)\,dt=
{(2a)^\mu(2b)^\nu\Gam(\mu+\nu+1)\ov(a^2+b^2)^{\mu+\nu+1}}.
\end{equation}
\end{lemma}

From this lemma, we obtain the following result
\begin{theorem}\label{th0}
If $|\mbox{Re}\ \nu|<1/2$ then
\begin{equation}\label{k5}
\hskip 0.75cm{\cal F}^{-1}[|\xi|^{\nu}K_{\nu}(|\xi|)](x)=
{2^{\nu-1}\Gam({n\ov2}+\nu)\ov\pi^{n\ov2}}\left(1+|x|^2\right)^{-{n\ov2}-\nu}.
\end{equation}
\end{theorem}

\ni{\bf Proof.} The integral in (\ref{k3}) is the same as the integral in (\ref{k4}), where $\mu=n/2-1,$
$a=|x|,$ and $b=1.$ Clearly the conditions of the lemma are satisfied and Theorem \ref{th0} follows at once.
$\Box$

We remark that, when $n=1,$ (\ref{k5}) becomes
$$
{\cal F}^{-1}[|\xi|^{\nu}K_{\nu}(|\xi|)](x)=
{2^{\nu-1}\Gam({1\ov2}+\nu)\ov\pi^{1\ov2}}\left(1+|x|^2\right)^{-{1\ov2}-\nu},
$$
a formula found in \cite{gs} and \cite{rg}. 

\bigskip
\ni{\bf 2. Inverse Fourier transforms of $|\xi|^\nu J_\nu(|\xi|)$ and $|\xi|^{-\nu} J_\nu(|\xi|).$}

As before, assume that $|\mbox{Re}\ \nu|<1/2.$ Formally, the inverse
Fourier transform of $|\xi|^{\nu}J_\nu(|\xi|)$ is
$$
{\cal F}^{-1}[|\xi|^{\nu}J_{\nu}(|\xi|)](x)=
  {1\ov(2\pi)^n}\int_{\mathbb R_n}e^{i\la x,\xi\ra}
  |\xi|^{\nu}J_\nu(|\xi|)\, d\xi.
$$
In general, the integral diverges at $\infty$ and so we introduce a
converging factor $e^{-\ep|\xi|}$ and take a limit as $\ep\to 0.$
Since $|\xi|^{\nu}J_\nu(|\xi|)$ is locally integrable, it defines a
tempered distribution and the limit exists in ${\cal S}'({\mathbb
R_n}),$ the space of tempered distributions on ${\mathbb R_n}.$ Thus
the precise meaning of the inverse Fourier transform is
\begin{equation}\label{f1}
{\cal F}^{-1}[|\xi|^{\nu}J_{\nu}(|\xi|)](x)=
  {1\ov(2\pi)^n}\lim_{\ep\to 0}\int_{\mathbb R_n}e^{i\la x,\xi\ra-\ep|\xi|}
  |\xi|^{\nu}J_\nu(|\xi|)\, d\xi.
\end{equation}
By introducing spherical coordinates, one can see that the integral on
the right-hand can be written as
$$
\int_0^\infty e^{-\ep r}r^{n+\nu-1}J_\nu(r)\{\int_{S_{n-1}}e^{i\la
rx,\om\ra}\, d\om\}\,dr.
$$
Since
\begin{equation}\label{f2}
\int_{S_{n-1}}e^{i\la rx,\om\ra}\, d\om=
  {(2\pi)^{\scs n\ov2}\ov|rx|^{{\scs n\ov2}-1}}J_{{\scs n\ov2}-1}(r|x|),
\end{equation}
it follows that
\begin{equation}\label{f3}
\hskip 0.75cm{\cal F}^{-1}[|\xi|^{\nu}J_{\nu}(|\xi|)](x)=
  {|x|^{1-{\scs n\ov2}}\ov(2\pi)^{\scs n\ov2}} \lim_{\ep\to
  0}\int_0^\infty e^{-\ep r}r^{{\scs n\ov2}+\nu}J_{{\scs
  n\ov2}-1}(r|x|)J_\nu(r)\, dr.
\end{equation}
In a similar manner we also have
\begin{equation}\label{f4}
\hskip 0.75cm{\cal F}^{-1}[|\xi|^{-\nu}J_{\nu}(|\xi|)](x)=
  {|x|^{1-{\scs n\ov2}}\ov(2\pi)^{\scs n\ov2}}\!  \lim_{\ep\to
  0}\int_0^\infty e^{-\ep r}r^{{\scs n\ov2}-\nu}J_{{\scs
  n\ov2}-1}(r|x|)J_\nu(r)\, dr.
\end{equation}

Both integrals appearing in formulas (\ref{f3}) and (\ref{f4}) are
particular cases of the following integral
\begin{equation}\label{f5}
I_\ep(a,b)=\int_0^\infty e^{-\ep t}t^{-\lam}J_\mu(at)J_\nu(bt)\, dt,
\end{equation}
where $a$ and $b$ are positive real numbers, $\lam,$ $\mu,$ and $\nu,$
complex numbers such that $\mbox{Re}(\mu+\nu+1)>\mbox{Re}(\lam).$ This
integral is a variant of the discontinuous integral of Weber and
Schafheitlin studied by Watson in his treatise on Bessel functions
\cite{wat}.

\begin{lemma}\label{l1}
As $\ep\to0$ $I_\ep(a,b)$ tends, in the sense of distributions, to
either one of the following limits:
\begin{equation}\label{f6}
\hskip .75 cm {b^\nu\Gam({\mu+\nu-\lam+1\ov2})\ov 2^\lam
  a^{\nu-\lam+1}\Gam(\nu+1)\Gam({\lam+\mu-\nu+1\ov2})}
  F({{\mu\!+\!\nu\!-\!\lam\!+\!1\ov2}},{{\nu\!-\!\lam\!-\!\mu\!+\!1\ov2}};\nu\!
  +\!1;{b^2\ov a^2}),
\end{equation}
if $0<b<a,$ or
\begin{equation}\label{f7}
\hskip .75 cm {a^\mu\Gam({\mu+\nu-\lam+1\ov2})\ov 2^\lam
  b^{\mu-\lam+1}\Gam(\mu+1)\Gam({\lam+\nu-\mu+1\ov2})}
  F({{\mu\!+\!\nu\!-\!\lam\!+\!1\ov2}},{{\mu\!-\!\lam\!-\!\nu\!+\!1\ov2}};
  \mu\!+\!1;{a^2\ov b^2}),
\end{equation}
if $0<a<b.$ 
\end{lemma}

\ni{\bf Proof.} We adapt to our situation the proof of the
Weber-Schafheitlin theorem as found in section 13.4 of Watson's
treatise \cite{wat}. It consists of expanding the integrand in
(\ref{f5}) in power series of $b,$ and passing to the limit as $\ep\to
0.$

1. Consider the case when $0<b<a.$ If we replace $b$ by $z,$ the
integral (\ref{f5}) is an analytic function of $z$ when $\mbox{Re}\,
z>0$ and $|\mbox{Im}\, z|<\ep.$ Introduce new constants $\alpha,$
$\beta,$ and $\gam$ defined by
$$
2\alpha=\mu+\nu-\lam+1,\qquad 2\beta=\nu-\lam-\mu+1,\qquad \gam=\nu+1
$$
or, equivalently,
$$
\lam=\gam-\alpha-\beta,\qquad \mu=\alpha-\beta,\qquad \nu=\gam-1,
$$
and rewrite (\ref{f5}) as follows:
\begin{eqnarray}\label{f5a}
\lefteqn{I_\ep(a,z) = \int_0^\infty e^{-\ep t}t^{\alpha+\beta-\gam}
J_{\alpha-\beta}(at)J_{\gam-1}(zt)\, dt } \\ & = &\int_0^\infty
e^{-\ep t}J_{\alpha-\beta}(at)
\Big\{\sum_{m=0}^\infty{(-1)^m(z/2)^{\gam+2m-1}t^{\alpha+\beta+2m-1}\ov
m!\Gam(\gam+m)}\Big\}\, dt \nonumber \\ & =
&\sum_{m=0}^\infty{(-1)^m(z/2)^{\gam+2m-1}\ov m!\Gam(\gam+m)}
\int_0^\infty e^{-\ep t}J_{\alpha-\beta}(at)t^{\alpha+\beta+2m-1}\,
dt. \nonumber
\end{eqnarray}
The interchange betwen the integration and summation signs is
justified because, when $|z|<\ep,$ the series
$$
\sum_{m=0}^\infty{(-1)^m(z/2)^{\gam+2m-1}\ov m!\Gam(\gam+m)}
\int_0^\infty e^{-\ep t}|J_{\alpha-\beta}(at)t^{\alpha+\beta+2m-1}|\, dt
$$
is absolutely convergent. 

2. We now evaluate the last integral on the right-hand side of
(\ref{f5a}). By expanding $J_{\alpha-\beta}(at)$ in power series and
integrating term by term we obtain
\begin{eqnarray*}
\lefteqn{\int_0^\infty e^{-\ep
 t}J_{\alpha-\beta}(at)t^{\alpha+\beta+2m-1}\, dt}\\ & &
 =\sum_{k=0}^\infty{(-1)^k(a/2)^{\alpha-\beta+2k}\ov
 k!\Gam(\alpha-\beta+k+1)} \int_0^\infty e^{-\ep
 t}t^{2\alpha+2(m+k)-1}\, dt \nonumber \\ & &
 =\sum_{k=0}^\infty{(-1)^k(a/2)^{\alpha-\beta+2k}\ov
 k!\Gam(\alpha-\beta+k+1)}
 {\Gam(2\alpha+2m+2k)\ov\ep^{2\alpha+2m+2k}}.\nonumber
\end{eqnarray*}
By using the duplication formula
$\Gam(2z)=2^{2z-1}\pi^{-1/2}\Gam(z)\Gam(z+{1\ov2})$ together with the
notation $(z)_k=z(z+1)\cdots(z+k-1)=\Gam(z+k)/\Gam(z)$ we rewrite the
last expression as:
\begin{eqnarray*}
\hspace{1cm}\lefteqn{\int_0^\infty e^{-\ep
  t}J_{\alpha-\beta}(at)t^{\alpha+\beta+2m-1}\, dt}\\ &= &
  {(a/2)^{\alpha-\beta}\ov(\ep^2)^{\alpha+m}}{\Gam(2\alpha+2m)\ov
    \Gam(\alpha-\beta+1)}
  \sum_{k=0}^\infty{(\alpha+m)_k(\alpha+m+{1\ov2})_k\ov
  k!(\alpha-\beta+1)_k} \left(-{a^2\ov c^2}\right)^k \nonumber \\ & =
  &
  {(a/2)^{\alpha-\beta}\ov(\ep^2)^{\alpha+m}}{\Gam(2\alpha+2m)\ov
\Gam(\alpha-\beta+1)}
  F(\alpha+m,\alpha+m+{1\ov2};\alpha-\beta+1;-{a^2\ov c^2}).\nonumber
\end{eqnarray*}
By using formula (\ref{hy5}) in the Appendix II:
$$
F(a,b;c;z)=(1-z)^{-a}F(a,c-b;c;{z\ov z-1}),
$$
we obtain that
\begin{eqnarray}\label{f5b}
\lefteqn{\hspace{2cm}\int_0^\infty e^{-\ep
 t}J_{\alpha-\beta}(at)t^{\alpha+\beta+2m-1}\, dt=}\\ & = &
 {(a/2)^{\alpha-\beta}\ov(a^2+\ep^2)^{\alpha+m}}{\Gam(2\alpha+2m)
  \ov\Gam(\alpha-\beta+1)}
 F(\alpha+m,\alpha+m+{1\ov2};\alpha-\beta+1;{a^2\ov
 a^2+\ep^2}).\nonumber
\end{eqnarray}
Substituting (\ref{f5b}) into (\ref{f5a}) we get
\begin{eqnarray}\label{f5c}
\lefteqn{\hspace{-4cm}I_\ep(a,z)=}
\end{eqnarray}
 $$
=\sum_{m=0}^\infty{(-1)^m(z/2)^{\gam+2m-1}(a/2)^{\alpha-\beta}\Gam(2\alpha+2m)
\ov m!\Gam(\gam+m)(a^2+\ep^2)^{\alpha+m}\Gam(\alpha-\beta+1)}
  F(\alpha+m,\alpha+m+{1\ov2};\alpha-\beta+1;{a^2\ov a^2+\ep^2}),
 $$
whenever $|z|<\ep .$ 

3. Following \cite{wat}, one can show that (\ref{f5c}) is valid
provided that $z$ satisfies the conditions
$$
\mbox{Re}(z)>0, \qquad |\mbox{Im}(z)|<\ep, \qquad
|z|<\sqrt{a^2+\ep^2}-\ep .
$$
Let $\delta >0$ be small enough so that
$b<\sqrt{a^2+\delta^2}-\delta,$ and take $0<\ep\leq\delta$ so that we
also have $b<\sqrt{a^2+\ep^2}-\ep$. In (\ref{f5c}) we may now let
$z=b$ and, when this is done, one can show, by the method of
majorants, that the resulting series converges uniformly with respect
to $\ep,$ $0<\ep\leq\delta,$ and therefore, as $\ep\to 0,$ the limit
of the series is equal to its value at $\ep=0.$ Thus
\begin{eqnarray}\label{f5d}
\lefteqn{\hspace{-4cm}\lim_{\ep\to0}\,I_\ep(a,b)=}
\end{eqnarray}
 $$
=\sum_{m=0}^\infty{(-1)^m(b/2)^{\gam+2m-1}(a/2)^{\alpha-\beta}\Gam(2\alpha+2m)
\ov m!\Gam(\gam+m)a^{2\alpha+2m}\Gam(\alpha-\beta+1)}
  F(\alpha+m,\alpha+m+{1\ov2};\alpha-\beta+1;1).
 $$
By formula (\ref{hy4}) in the Appendix II, we have
$$
F(\alpha+m,{1\ov2}-\beta-m;\alpha-\beta+1;1)=
   {\Gam(\alpha-\beta+1)\Gam(1/2)\ov\Gam(1-\beta-m)\Gam(\alpha+m+{1\ov2})}.
$$
On the other hand,
$$
\Gam(1/2)\Gam(2\alpha+2m)=2^{2\alpha+2m-1}\Gam(\alpha+m)\Gam(\alpha+m+{1\ov2}),
$$
$$
\Gam(1-\beta-m)={(-1)^m\pi\csc(\pi\beta)\ov\Gam(\beta+m)},
$$
$$
\Gam(1-\beta)\Gam(\beta)=\pi\csc(\pi\beta),
$$
therefore
$$
{\Gam(2\alpha+2m)\Gam(1/2)\ov
   \Gam(\gam+m)\Gam(1-\beta-m)\Gam(\alpha+m+{1\ov2})} =
{2^{2\alpha+2m-1}\Gam(\alpha+m)\Gam(\beta+m)\ov(-1)^m\Gam(1-\beta)\Gam(\beta)
 \Gam(\gam+m)}.
$$
Substituting these formulas into (\ref{f5d}), we get
\begin{eqnarray}\label{f5e}
\lefteqn{\lim_{\ep\to0}\,I_\ep(a,b)}\\
&=& \sum_{m=0}^\infty{b^{\gam-1}\Gam(\alpha+m)\Gam(\beta+m)\ov
 2^{\gam-\alpha-\beta}a^{\alpha+\beta}m!\Gam(1-\beta)\Gam(\beta)\Gam(\gam+m)}
\left(b^2\ov a^2\right)^m\nonumber\\
&=& {b^{\gam-1}\Gam(\alpha)\ov
 2^{\gam-\alpha-\beta}a^{\alpha+\beta}\Gam(\gam)\Gam(1-\beta)}
F(\alpha,\beta;\gam;{b^2\ov a^2}).\nonumber
\end{eqnarray}
Finally, returning to the constants $\mu,$ $\nu,$ and $\lam,$ we
obtain the expression (\ref{f6}) in the first part of the lemma.

4. In the case $0<a<b,$ we proceed in an analogous manner to obtain
the expression (\ref{f7}) in the second part of the lemma. $\Box$

\medskip
We now use Lemma \ref{l1} to evaluate the inverse Fourier transforms
of $|\xi|^\nu J_\nu(|\xi|)$ and $|\xi|^{-\nu} J_\nu(|\xi|)$
respectively defined by formulas (\ref{f3}) and (\ref{f4}).

\begin{theorem}\label{th3} The inverse Fourier transform of 
$|\xi|^{\nu}J_\nu(|\xi|)$ is the distribution defined by
\begin{eqnarray}\label{f5f}
\hspace{1.5cm}\lefteqn{{\cal F}^{-1}[|\xi|^{\nu}J_\nu(|\xi|)](x)
=}\hspace{1.5cm} \\ & & =\left\{\begin{array}{ll} \ds\ \
\sin({n\pi\ov2}){2^\nu\Gamma({n\ov2}+\nu)\ov\pi^{{n\ov2}+1}}
\,(1-|x|^2)^{-{n\ov2}-\nu}, & 0<|x|<1,\\ \\
\ds-\sin(\nu\pi){2^\nu\Gamma({n\ov2}+\nu)\ov\pi^{{n\ov2}+1}}\,
(|x|^2-1)^{-{n\ov2}-\nu}, & 1<|x|.
\end{array} \right.\nonumber
\end{eqnarray}
\end{theorem}

\ni{\bf Proof.} The integral on the right-hand side of (\ref{f3})
corresponds to the integral in (\ref{f5}) where
$$
\mu={n\ov2}-1,\quad \nu=\nu,\quad  \lam=-{n\ov2}-\nu, \quad a=|x|, \quad b=1.
$$
If $1<|x|,$ then from formula (\ref{f6}) we obtain for the limit in
 (\ref{f3}): \begin{eqnarray*} 
\lefteqn{\lim_{\ep\to0}\int_0^\infty
 e^{-\ep r}r^{{n\ov2}+\nu}J_{{n\ov2}-1}(r|x|)J_\nu(r)\,dr =}\\ & &
 ={\Gam({n\ov2}+\nu)\ov2^{-{n\ov2}-\nu}|x|^{{n\ov2}+2\nu+1}\Gam(\nu+1)
     \Gam(-\nu)}
 F({n\ov2}+\nu,\nu+1;\nu+1;{1\ov|x|^2}).
\end{eqnarray*}
From the known relation for hypergeometric series
\begin{equation}\label{f5g}
F(a,b;c;z)=(1-z)^{c-a-b}F(c-a,c-b;c;z)
\end{equation}
we have that
$$
F({n\ov2}+\nu,\nu+1;\nu+1;{1\ov |x|^2})=\left(|x|^2-1\ov
|x|^2\right)^{-{n\ov2}-\nu}.
$$
On the other hand
$$
\Gam(\nu+1)\Gam(-\nu)={-\pi\ov\sin(\pi\nu)}.
$$
Thus, for $1<|x|,$ the above limit is equal to
\begin{eqnarray}\label{f5i}
\lefteqn{\lim_{\ep\to0}\int_0^\infty e^{-\ep r}\!r^{{n\ov2}+\nu}
  \!J_{{n\ov2}-1}(r|x|)J_\nu(r)\,dr =}\\
& & = -\sin(\pi\nu){|x|^{{n\ov2}-1}\Gam({n\ov2}+\nu)\ov2^{-{n\ov2}-\nu}\pi}
   \left(|x|^2-1\right)^{-{n\ov2}-\nu}, \nonumber
\end{eqnarray}

If $0<|x|<1,$ then formula (\ref{f7}) yields
 \begin{eqnarray}\label{f5j} 
\lefteqn{\lim_{\ep\to0}\int_0^\infty
 e^{-\ep r}r^{{n\ov2}+\nu}J_{{n\ov2}-1}(r|x|)J_\nu(r)\,dr =}\\ & &
 ={|x|^{{n\ov2}-1}\Gam({n\ov2}+\nu)\ov2^{-{n\ov2}-\nu}\Gam({n\ov2})
    \Gam(1-{n\ov2})}
 F({n\ov2}+\nu,{n\ov2};{n\ov2};|x|^2) \nonumber\\ &
 &=\sin({n\pi\ov2}){|x|^{{n\ov2}-1}\Gam({n\ov2}+\nu)\ov2^{-{n\ov2}-\nu}\pi}
 \left(1-|x|^2\right)^{-{n\ov2}-\nu}, \nonumber
\end{eqnarray}
after making the replacements
$$
\Gam({n\ov2})\Gam(1-{n\ov2})={\pi\ov\sin({n\pi/2})}
$$
and 
$$
F({n\ov2}+\nu,{n\ov2};{n\ov2};|x|^2)=(1-|x|^2)^{-{n\ov2}-\nu}.
$$
The last expression results from (\ref{f5g}). Substituting (\ref{f5i})
and (\ref{f5j}) into (\ref{f3}) we obtain (\ref{f5f}) and the theorem
is proved $\Box$

By reasoning in the same manner the following result, whose proof is
left to the reader, holds.
\begin{theorem}\label{th4} The inverse Fourier transform of 
$|\xi|^{-\nu}J_\nu(|\xi|)$ is the distribution defined by
\begin{equation}\label{f5h}
{\cal F}^{-1}[|\xi|^{-\nu}J_\nu(|\xi|)](x) = 
\left\{\begin{array}{ll}
\ds{(1-|x|^2)^{\nu-{n\ov2}}\ov2^\nu\pi^{n\ov2}\Gamma(\nu-{n\ov2}+1)},
       & 0<|x|<1,\\
\\
\hspace{1.5cm} 0, 
       & 1<|x|. 
\end{array} \right.\nonumber
\end{equation}
\end{theorem}

\medskip
{\bf Remark.} From theorem \ref{th4} we immediately derive a formula
for the inverse Fourier transform of $|\xi|^\nu J_{-\nu}(|\xi|),$
namely,
\begin{equation}\label{f5h1}
{\cal F}^{-1}[|\xi|^{\nu}J_{-\nu}(|\xi|)](x) = 
\left\{\begin{array}{ll}
\ds{(1-|x|^2)^{-\nu-{n\ov2}}\ov2^{-\nu}\pi^{n\ov2}\Gamma(1-\nu-{n\ov2})},
       & 0<|x|<1,\\
\\
\hspace{1.5cm} 0, 
       & 1<|x|. 
\end{array} \right.\nonumber
\end{equation}

Since 
$$
\Gam(1-\nu-{n\ov2})\Gam(\nu+{n\ov2})={\pi\ov\sin[(\nu+{n\ov2})\pi]},
$$
we may rewrite (\ref{f5h1}) as
\begin{eqnarray}\label{f5h2}
\hspace{2cm}\lefteqn{{\cal F}^{-1}[|\xi|^{\nu}J_{-\nu}(|\xi|)](x) =
}\\ &=&\left\{\begin{array}{ll}
\ds{2^\nu\sin[(\nu+{n\ov2})\pi]\Gam(\nu+{n\ov2})\ov\pi^{{n\ov2}+1}}
  (1-|x|^2)^{-\nu-{n\ov2}},
& 0<|x|<1,\\ \\
\hspace{1.5cm} 0, 
       & 1<|x|. 
\end{array} \right.\nonumber
\end{eqnarray}

\medskip
\ni{\bf 3. The inverse Fourier transform of $|\xi|^\nu N_\nu(|\xi|).$}

Theorems \ref{th3} and \ref{th4} imply the following result
\begin{theorem}\label{th5} The inverse Fourier transform of 
$|\xi|^\nu N_\nu(|\xi|)$ is the distribution defined by
\begin{eqnarray}\label{f5k}
\hspace{1.5cm}\lefteqn{{\cal F}^{-1}[|\xi|^{\nu}N_\nu(|\xi|)](x) =}
  \hspace{1.5cm} \\
& & =\left\{\begin{array}{ll}
\ds -\cos({n\pi\ov2}){2^\nu\Gamma({n\ov2}+\nu)\ov\pi^{{n\ov2}+1}}
      \,(1-|x|^2)^{-{n\ov2}-\nu},
       & 0<|x|<1,\\
\\
 \ds-\cos(\nu\pi){2^\nu\Gamma({n\ov2}+\nu)\ov\pi^{{n\ov2}+1}}\,
(|x|^2-1)^{-{n\ov2}-\nu}, 
       & 1<|x|. 
\end{array} \right.\nonumber
\end{eqnarray}
\end{theorem}
\ni{\bf Proof.} Recall that
$$
N_\nu(|\xi|)={J_\nu(|\xi|)\cos(\nu\pi)-J_{-\nu}(|\xi|)\ov\sin(\nu\pi)}.
$$
Thus, from formulas (\ref{f5f}) and (\ref{f5h2}) we get
\begin{eqnarray*}
\hspace{1.5cm}\lefteqn{{\cal F}^{-1}[|\xi|^{\nu}N_\nu(|\xi|)](x) =}
  \hspace{1.5cm} \\
& & =\left\{\begin{array}{ll}
\ds {\cos(\nu\pi)\ov\sin(\nu\pi)}{2^\nu\Gamma({n\ov2}+\nu)\ov\pi^{{n\ov2}+1}}
      \,(1-|x|^2)^{-{n\ov2}-\nu},
       & 0<|x|<1,\\
\\
 \ds-\cos(\nu\pi){2^\nu\Gamma({n\ov2}+\nu)\ov\pi^{{n\ov2}+1}}\,
(|x|^2-1)^{-{n\ov2}-\nu}, 
       & 1<|x|, 
\end{array} \right. \\
\\
& & -\left\{\begin{array}{ll}
\ds
{\sin[(\nu+{n\ov2})\pi]\ov\sin(\nu\pi)}{2^\nu\Gamma({n\ov2}+\nu)\ov
  \pi^{{n\ov2}+1}}
      \,(1-|x|^2)^{-{n\ov2}-\nu},
       & 0<|x|<1,\\
\\
 \hspace{3cm} 0, & 1<|x|. 
\end{array} \right.  
\end{eqnarray*}
Now
$$
\sin[(\nu+{\pi\ov2})]\pi=\sin(\nu\pi)\cos({n\pi\ov2})+
  \cos(\nu\pi)\sin({n\pi\ov2})
$$
and so formula (\ref{f5k}) follows at once. $\Box$

\section{A fundamental solution with support in $\ovl{D_-^n}$}
  \setcounter{equation}{0}\label{hy} 
We return to the problem (\ref{eq3}) and, in order to obtain a formula
for a fundamental solution with support in the hyperbolic region, we
have to modify formula (\ref{eq8}). Consider the function
\begin{equation}\label{eq10}
\td F_-(\xi,y) = \left\{\begin{array}{ll}
\ds{3^{2/3}\Gamma(4/3)\over2^{1/3}}({t\over|\xi|})^{1/3}J_{1/3}(t|\xi|)
       & \mbox{for $y\leq 0$}\\
\\
\ds \hskip 1cm 0 & \mbox{for  $y\geq 0$}, 
\end{array} \right.
\end{equation}
where $t=(2/3)(-y)^{3/2}.$ From formula (\ref{b7}) in the Appendix I,
it follows that the limit of $\td F_-(\xi,y)$ as $y\to0-$ is equal to
zero and so condition (i) in Section \ref{pr} holds. Also it follows
from the same formula that
$$
\partial_y\{({t\over|\xi|})^{1/3}J_{1/3}(t|\xi|)\}_{|y=0}=
-{2^{1/3}\over3^{2/3}\Gamma(4/3)},
$$
hence, the $y$-derivative of $\td F_-(\xi,y)$ at $y=0-$ is equal to
$-1,$ and so condition (ii) is also satisfied. By calculating the
inverse Fourier transform of $\td F_-(\xi,y)$ we have the following
result:
\begin{theorem}\label{th1}The inverse Fourier transform of 
$\td F_-(\xi,y)$ is the distribution
\begin{equation}\label{eq11}
F_-(x,y) = \left\{\begin{array}{ll}
\ds{3^n\Gamma({4\ov3})\ov2^{2/3}\pi^{n/2}\Gamma({4\ov3}-{n\ov2})}
      \,|9|x|^2+4y^3|^{{1\ov3}-{n\ov2}}
       & \mbox{in $D_-^n$}\\
\\
 \hskip 1cm 0 & \mbox{elsewhere,} 
\end{array} \right.
\end{equation}
where $D_-^n=\{(x,y)\in{\mathbb R}^{n+1}: 9|x|^2+4y^3<0\}.$ It is a
fundamental solution for the operator (\ref{eq1}) whose support is the
closure of the region $D_-^n.$
\end{theorem}
 
\ni{\bf Proof.} We must evaluate the inverse Fourier transform of
$(t/|\xi|)^{1/3}J_{1/3}(t|\xi|).$ Recall that if $G(x)$ is the inverse
Fourier transform of $f(\xi),$ then $(1/a^n)G(x/a)$ is the inverse
Fourier transform of $f(a\xi).$ By appying formula (\ref{f5h}) for
$\nu=1/3,$ we obtain
\begin{eqnarray*}
\hspace{1cm}\lefteqn{{\cal F}^{-1}[({t/|\xi|})^{1/3}J_{1/3}(t|\xi|)](x) =}
\hspace{1.5cm}\\ 
& &=\left\{\begin{array}{ll}
\!\ds{t^{{2\ov3}-n}\ov2^{1/3}\pi^{n\ov2}\Gamma({4\ov3}-{n\ov2})}\!
\left(1-{|x|^2\ov t^2}\right)^{{1\ov3}-{n\ov2}}\!\!, &\! 0<|x|<t,\\ 
\\
\hspace{1.5cm} 0, 
       & t<|x|. \\
\end{array} \right.\\
\\
& & = \left\{\begin{array}{ll}
\!\ds{1\ov2^{1/3}\pi^{n\ov2}\Gamma({4\ov3}-{n\ov2})}\!
  \left(t^2-|x|^2\right)^{{1\ov3}-{n\ov2}}\!\!,
       &\! 0<|x|<t,\\
\\
\hspace{1.5cm} 0, 
       & t<|x|.
\end{array} \right.\\
\\
& & =\left\{\begin{array}{ll}
\!\ds{3^{n-{2\ov3}}\ov2^{1/3}\pi^{n\ov2}\Gamma({4\ov3}-{n\ov2})}\!
   \left|9|x|^2+4y^3\right|^{{1\ov3}-{n\ov2}}
       & \mbox{in}\ D_-^n\\
\\
\hspace{1.5cm} 0, 
       & \mbox{elsewhere}.
\end{array} \right.
\end{eqnarray*}

Multiplication by the constant $3^{2/3}\Gam(4/3)/2^{1/3}$ in formula
(\ref{eq10}) yields (\ref{eq11}) which proves the theorem. $\Box$

\medskip
{\bf Remarks 1.} If $n=1,$ then $F_-(x,y)$ coincides with the
distribution defined by formula (\ref{eq2b}).

{\bf 2.} For all values of $n$ the support of $F_-(x,y)$ is the
closure of the region $ D_-^n.$ This follows from results in \cite{gs}
about the generalized function $P^\lambda_+$ where $P(x)$ is the
quadratic polynomial
$$
P(x)=x_1^2+\ldots+x_p^2-x_{p+1}^2-\ldots -x_{p+q}^2,
$$
$p, q\geq 1.$ See also the example in \cite{sch} about the
distribution $\mbox{Pf}.\, s^\lambda,$ where $s$ denotes the
hyperbolic distance in $\mathbb{R}^N,$ $N>1.$

\section{A fundamental solution with support in $\ovl{D_+^n}$}
  \setcounter{equation}{0}\label{el}
Consider the function $\td F_+(\xi,y)$ defined by
\begin{equation}\label{eq12}
\td F_+(\xi,y) = \left\{\begin{array}{ll}
\ds \gamma\cdot({s\over|\xi|})^{1/3}K_{1/3}(s|\xi|)
    & \mbox{if $y\geq 0$} \\
\ds
\delta\cdot({t\over|\xi|})^{1/3}N_{-1/3}(t|\xi|)
       & \mbox{if $y\leq 0$},
\end{array} \right.
\end{equation}
where $s=(2/3)y^{3/2},$ $t=(2/3)(-y)^{3/2},$ and the constants
$\gamma$ and $\delta$ are respectively given by
\begin{equation}\label{eq12a}
\gamma=-{2^{2/3}\ov3^{4/3}\Gamma(2/3)}
\qquad\mbox{and}\qquad
\delta={2\pi\over2^{1/3}3^{4/3}\Gamma(2/3)},
\end{equation}
and where $N_{-1/3}$ is the Neumann function defined by (\ref{n1}). By
using the formulas in the Appendix I, it is a matter of verification
that the conditions (i) and (ii) in Section \ref{pr} are
satisfied. Thus its inverse Fourier transform, denoted by
$F^\sharp(x,y),$ defines a fundamental solution for the operator
(\ref{eq1}). The following theorem gives an explicit expression for
$F^\sharp(x,y).$
\begin{theorem}\label{th2}The inverse Fourier transform of 
$\td F_+(\xi,y)$ is the distribution defined by
\begin{equation}\label{eq13}
\hskip 0.75cm F^\sharp(x,y) = \left\{\begin{array}{ll}
\ds-{3^{n-2}\ov
2^{2/3}\pi^{n/2}}{\Gamma({n\ov2}-{1\ov3})\ov\Gamma(2/3)}
      \,(9|x|^2+4y^3)^{{1\ov3}-{n\ov2}}
       & \mbox{in $D_+^n$}\\
\\
 \ds-\cos({n\pi\ov2}){2^{1/3}3^{n-2}\ov\pi^{n/2}}
  {\Gamma({n\ov2}-{1\ov3})\ov\Gamma(2/3)}
    \, |9|x|^2+4y^3|^{{1\ov3}-{n\ov2}}
     & \mbox{in $D_-^n$.} 
\end{array} \right.
\end{equation}
\end{theorem}

\ni{\bf Proof.} 1. We start by evaluating the inverse Fourier
transform of $(s/|\xi|)^{1/3}K_{1/3}(s|\xi|),$ $s=2y^{3/2}/3,$
$y\geq0.$ Since $K_\nu(z)=K_{-\nu}(z),$ it follows from formula
(\ref{k5}) that
$$
{\cal F}^{-1}[(s/|\xi|)^{1/3}K_{1/3}(s|\xi|)] =
    {2^{-{4\ov3}}\Gam({n\ov2}-{1\ov3})\ov\pi^{n/2}}
\left(s^2+|x|^2\right)^{{1\ov3}-{n\ov2}}.
$$
After multiplying by the constant $\gam$ in (\ref{eq12a}) and
reintroducing the variables $x$ and $y,$ we can see that the
right-hand side of the last expression is equal to
\begin{equation}\label{eq13a}
-{3^{n-2}\Gam({n\ov2}-{1\ov3})\ov2^{2/3}\pi^{n/2}\Gam(2/3)}
 \left(9|x|^2+4y^3\right)^{{1\ov3}-{n\ov2}},
\quad y\geq0,
\end{equation}
that is $F^\sharp(x,y)$ in $D_+^n\cap\{y\geq0\}.$

2. Next we evaluate the inverse Fourier transform of 
$(t/|\xi|)^{1/3}N_{-1/3}(t|\xi|),$ $t=2(-y)^{3/2}/3,$ $y\leq0.$ From
formula (\ref{f5k}) we obtain
\begin{eqnarray*}
\hspace{1.5cm}\lefteqn{{\cal F}^{-1}[(t|\xi|)^{-1/3}N_{-1/3}(t|\xi|)]=}
\hspace{1.5cm} \\
& & =\left\{\begin{array}{ll}
\ds -\cos({n\pi\ov2}){2^{-1/3}\Gamma({n\ov2}-{1\ov3})\ov\pi^{{n\ov2}+1}}
      \,(t^2-|x|^2)^{{1\ov3}-{n\ov2}},  & 0<|x|<t,\\
\\
 \ds-\left(1\ov2\right){2^{-1/3}\Gamma({n\ov2}-{1\ov3})\ov\pi^{{n\ov2}+1}}\,
   (|x|^2-t^2)^{{1\ov3}-{n\ov2}},  & t<|x|. 
\end{array} \right.
\end{eqnarray*}
After multiplying both sides by the constant $\delta$ in (\ref{eq12a})
and reverting to the variables $x$ and $y,$ we can see that the
right-hand side of the last expression can be written as
\begin{equation}\label{eq13b}
\hspace{.75cm}\left\{\begin{array}{ll}
\ds -{3^{n-2}\Gamma({n\ov2}-{1\ov3})\ov2^{2/3}\pi^{n/2}\Gam(2/3)}
      \,\left(9|x|^2+4y^3\right)^{{1\ov3}-{n\ov2}},  
    & \mbox{in}\ D_+^n\cap\{y\leq0\}\\ \nonumber
\\
 \ds-\cos({n\pi\ov2}){2^{1/3}3^{n-2}\Gamma({n\ov2}-{1\ov3})\ov
  \pi^{n/2}\Gam(2/3)}\,
   \left|9|x|^2+4y^3\right|^{{1\ov3}-{n\ov2}},  & \mbox{in}\ D_-^n, 
\end{array} \right.
\end{equation}
and so (\ref{eq13b}) coincides with $F^\sharp(x,y)$ in
$D_+^n\cap\{y\leq0\}\cup D_-^n.$ Therefore, from (\ref{eq13a}) and
(\ref{eq13b}) we obtain (\ref{eq13}) and the theorem is proved. $\Box$

\medskip
{\bf Remarks 1.} In view of the results in \cite{gs} relative to the
generalized function $P_+^\lambda$ (see Remark 2 after the proof of
Theorem \ref{th1}) the supports of the distributions
$$
(9|x|^2+4y^3)^{{1\ov3}-{n\ov2}}\qquad\mbox{and}\qquad 
  |9|x|^2+4y^3|^{{1\ov3}-{n\ov2}}
$$
that appear in formula (\ref{eq13}) are the closures of $D_+^n$ and
$D_-^n,$ respectively. In other words, for no value of $n$ can these
supports be just the boundaries of these regions.

{\bf 2.} If $n$ is {\em odd}, then $\cos(n\pi/2)=0,$ and we rewrite
formula (\ref{eq13}), using the notation $F_+(x,y)$ instead of 
$F^\sharp(x,y),$ as follows: 
\begin{equation}\label{eq14}
\hskip 0.75cm F_+(x,y) = \left\{\begin{array}{ll}
\ds-{3^{n-2}\ov
2^{2/3}\pi^{n/2}}{\Gamma({n\ov2}-{1\ov3})\ov\Gamma(2/3)}
      \,(9|x|^2+4y^3)^{{1\ov3}-{n\ov2}}
       & \mbox{in $D_+^n$}\\
\\
   \hskip 1cm 0  &\mbox{elsewhere.} 
\end{array} \right.
\end{equation}
This is a fundamental solution whose support is $\ovl{D_+^n}.$ In
particular, if $n=1$ we obtain formula (\ref{eq2a}).

{\bf 3.} If $n$ is {\em even}, then $F^\sharp(x,y)$ is not necessarily
identically zero in $D_-^n$ and its support may be the whole of
${\mathbb R}^{n+1}.$ Suppose that $n=2k,\ k>0,$ and let us compare the
constant in formula (\ref{eq13}), relative to the region $D_-^n,$ to
the constant in the expression (\ref{eq11}) of $F_-(x,y).$ Let
$$
A={(-1)^{k+1}2^{1/3}3^{2(k-1)}\ov\pi^k}{\Gamma(k-{1\ov3})\ov\Gamma(2/3)}
$$
be the constant in (\ref{eq13}) and let
$$
B={3^{2k}\ov2^{2/3}\pi^k}{\Gamma(4/3)\ov\Gamma({4\ov3}-k)}=
    {(-1)^{k+1}3^{2k-1}\ov2^{2/3}\pi^k}{\Gamma(k-{1\ov3})\ov\Gamma(2/3)}
$$
be the constant in (\ref{eq11}). Since $3A-2B=0$ it follows that the
distribution
\begin{equation}\label{eq15}
F_+(x,y)=3F^\sharp(x,y)-2F_-(x,y)
\end{equation}
is now a fundamental solution for the operator (\ref{eq1}) supported
in $\ovl{D_+^n}$ and we have for this $F_+(x,y)$ the same expression
as that of (\ref{eq14}). In conclusion, for all values on $n$ we
always get two fundamental solutions: one whose support is the closure
of $D_+^n$ and another whose support is the closure of $D_-^n.$

\section{Appendix}\setcounter{equation}{0}\label{ap}
{\bf I. Bessel functions}

The function $J_\nu(z)$ of a complex variable $z$ defined by
\begin{equation}\label{b1}
J_\nu(z)=\sum_{r=0}^\infty{(-1)^rz^{\nu+2r}
 \over2^{\nu+2r}r!\Gamma(\nu+r+1)},\quad |z|<\infty,\quad |\arg z|<\pi, 
\end{equation}
is called the {\em Bessel function of the first kind of order $\nu$}.

We also need the Bessel functions $I_\nu(z)$ and $K_\nu(z)$ defined by
\begin{equation}\label{b2}
I_\nu(z)=\sum_{r=0}^\infty{z^{\nu+2r}
 \over2^{\nu+2r}r!\Gamma(\nu+r+1)},\quad |z|<\infty,\quad |\arg z|<\pi,
\end{equation}
and
\begin{equation}\label{b3}
K_\nu(z)={\pi\csc(\nu\pi)\over2}\{I_{-\nu}(z)-I_\nu(z)\},\quad
    \nu\neq 0,\pm1, \pm2,\dots
\end{equation}
as well as the {\em Neumann function} 
\begin{equation}\label{n1}
N_\nu(z)={J_\nu(z)\cos(\nu\pi)-J_{-\nu}(z)\ov\sin(\nu\pi)}.
\end{equation}

Note that throughout this work, we only deal with Bessel functions of
order $\pm1/3$.  Recall that $Ai(z)$ was defined in formula
(\ref{eq5}) by
$$
Ai(z)={z^{1/2}\over 3}[I_{-1/3}({2\over3}z^{3/2})-I_{1/3}({2\over3}z^{3/2})]
     = {1\over\pi}({z\over 3})^{1/2}K_{1/3}({2\over3}z^{3/2}).
$$
If we set $s={2\over3}z^{3/2}$, then we may rewrite $Ai(z)$ as 
\begin{equation}\label{b4}
Ai(z)={1\over3^{2/3}2^{1/3}}s^{1/3}[I_{-1/3}(s)-I_{1/3}(s)].
\end{equation}
From the series expansion of $I_\nu(z)$ it follows that
\begin{equation}\label{b4a}
s^{1/3}I_{-1/3}(s)=
  {1\over2^{-1/3}\Gamma(2/3)} + {s^2\over2^{5/3}\Gamma(5/3)}+\ldots
\end{equation}
and
\begin{equation}\label{b4b}
s^{1/3}I_{1/3}(s)=
  {s^{2/3}\over2^{1/3}\Gamma(4/3)} + {s^{8/3}\over2^{7/3}\Gamma(7/3)}+\ldots
\end{equation}
Consequently from (\ref{b4}), (\ref{b4a}), and (\ref{b4b}) we obtain
$$
Ai(0)={3^{-2/3}\over\Gamma(2/3)},
$$
the first expression in formula (\ref{eq6a}). Similarly, by
differentiating $Ai(z)$ and setting $z=0$, we get the second
expression in (\ref{eq6a})
$$
Ai'(0)=-{3^{-4/3}\over\Gamma(4/3)}.
$$

In an analogous way, recall that $Bi(z)$ in (\ref{eq5}) was defined by
$$
Bi(z)=({z\over3})^{1/2}[I_{-1/3}({2\over3}z^{3/2})+I_{1/3}({2\over3}z^{3/2})].
$$
If we set $s={2\over3}z^{3/2}$, then the last expression becomes
\begin{equation}\label{b5}
Bi(z)={1\over2^{1/3}3^{1/6}}s^{1/3}[I_{-1/3}(s)+I_{1/3}(s)],
\end{equation}
and again from (\ref{b4a}) and (\ref{b4b}) we obtain
$$
Bi(0)={3^{-1/6}\over\Gamma(2/3)},\qquad Bi'(0)={3^{-5/6}\over\Gamma(4/3)},
$$
which are the two expressions in (\ref{eq6b}).

Finally, from these results we obtain the value of the Wronskian of $Ai(z)$
and $Bi(z)$ at $z=0$:
$$
W(Ai(z),Bi(z))_{|_{z=0}}={2\cdot3^{-3/2}\over\Gamma(2/3)\Gamma(4/3)}=1/\pi,
$$
because $\Gamma(2/3)\Gamma(4/3)=2\pi\cdot3^{-3/2}.$

For future reference, we also need the following expansions derived from the
series that defines $J_\nu(z)$:
\begin{equation}\label{b6}
t^{1/3}J_{-1/3}(t)=
  {1\over2^{-1/3}\Gamma(2/3)}-{t^2\over2^{5/3}\Gamma(5/3)}+\ldots
\end{equation}
and
\begin{equation}\label{b7}
t^{1/3}J_{1/3}(t)=
  {t^{2/3}\over2^{1/3}\Gamma(4/3)}-{t^{8/3}\over2^{7/3}\Gamma(7/3)}+\ldots
\end{equation}
\medskip

\ni{\bf II. Hypergeometric series and functions}

Let $a,$ $b,$ and $c$ be arbitrary complex numbers and let $z$ be a
complex variable. The power series
\begin{equation}\label{hy1}
F(a,b;c;z)=\sum_{n=0}^\infty{(a,n)(b,n)\ov(c,n)}{z^n\ov n!},
\end{equation}
where 
\begin{equation}\label{hy1a}
\hskip 0.75cm (a,0)=1,\
(a,n)={\Gam(a+n)\ov\Gam(a)}=a(a+1)\cdots(a+n-1),\ n=1,2,\ldots
\end{equation}
and we assume that $c\neq 0,-1,-2,\ldots,$ is called a {\em
hypergeometric series}. It is known \cite{ww} that the series
(\ref{hy1}) is a solution, valid near $z=0,$ of the {\em
hypergeometric equation}
\begin{equation}\label{hy2}
z(1-z){d^2u\ov dz^2}+\{c-(a+b+1)z\}{du\ov dz}-ab\,u=0,
\end{equation}
for which every point is an ordinary point, except $0,$ $1,$ and
$\infty,$ that are regular singular points.

If either $a$ or $b$ is a negative integer, then the series
(\ref{hy1}) terminates; if $c$ is a negative integer, the series is
meaningless because all terms after the $(1-c)$\,th have a zero
denominator. As it is known \cite{erd}, it is possible to redefine the
series (\ref{hy1}) so that it still is a solution of the
hypergeometric equation. We exclude this possibility from our
considerations because, in the cases that interest us, $c$ is never a
negative integer.

The hypergeometric series is absolutely convergent for $|z|<1$ and so
defines, in the open disk, an analytic function of $z$ which is
regular at $z=0.$ The point $z=1$ is however a branch point and if a
cut is made from $1$ to $+\infty$ along the $x$-axis, it can be shown
\cite{ww} that series can be continued analytically and defines an
analytic function throughout the cut plane that we still denote by
$F(x,b;c;z).$ If $\mbox{Re}(c)>\mbox{Re}(b)>0,$ this analytic
extension can be represented by Euler's formula
\begin{equation}\label{hy3}
F(a,b;c;z)={\Gam(c)\ov\Gam(b)\Gam(c-b)}
\int_0^1t^{b-1}(1-b)^{c-b-1}(1-tz)^{-a}\, dt,
\end{equation}
for $|\mbox\,{arg}(1-z)|<\pi.$

In general, the hypergeometric series (\ref{hy1}) diverges for
$|z|=1.$ However, if $\mbox{Re}(c-a-b)>0,$ we have absolute
convergence for $|z|=1.$ Moreover,
\begin{equation}\label{hy4}
F(a,b;c;1)={\Gam(c)\Gam(c-a-b)\ov\Gam(c-a)\Gam(c-b)}.
\end{equation}

Hypergeometric functions satisfy among themselves quite a number of
important relations of which we just list the following two that we
used in Section \ref{ft}:
\begin{equation}\label{hy5}
F(a,b;c;z)=(1-z)^{-a}F(a,c-b;c;{z\ov z-1})
\end{equation}
and
\begin{equation}\label{hy6}
F(a,b;c;z)=(1-z)^{c-a-b}F(c-a,c-b;c;z).
\end{equation}
For a complete list of such relations, the reader should consult
Erd\'ely \cite{erd}.

\end{document}